\documentclass[]{./layout/interact}
\usepackage{./layout/Layout}
\graphicspath{./{figures/}}

\title{The Buildings Gallery: visualising buildings}

\author{
\name{Bram Bekker\textsuperscript{a,b} and Maarten Solleveld\textsuperscript{a,c}}
\affil{\textsuperscript{a}Radboud Universiteit,
Department of Mathematics,
Heyendaalseweg 135,
6525AJ Nijmegen,
The Netherlands\newline
\textsuperscript{b}bekker.math@gmail.com\newline
\textsuperscript{c}m.solleveld@science.ru.nl}}
\date{Last modified: \today}

\begin{document}

\maketitle

\begin{abstract}
Buildings are beautiful mathematical objects tying a variety of subjects in algebra and geometry together in a very direct sense. They form a natural bridge to visualising more complex principles in group theory. As such they provide an opportunity to talk about the inner workings of mathematics to a broader audience, but the visualisations could also serve as a didactic tool in teaching group and building theory, and we believe they can even inspire future research.

We present an algorithmic method to visualise these geometric objects. The main accomplishment is the use of existing theory to produce three dimensional, interactive models of buildings associated to groups with a $BN$-pair. The final product, an interactive web application called The Buildings Gallery, can be found at  https://buildings.gallery \citep{Bekker2020}.
\end{abstract}


\begin{figure}[h!]
    \centering
    \includegraphics[width = 0.95\linewidth]{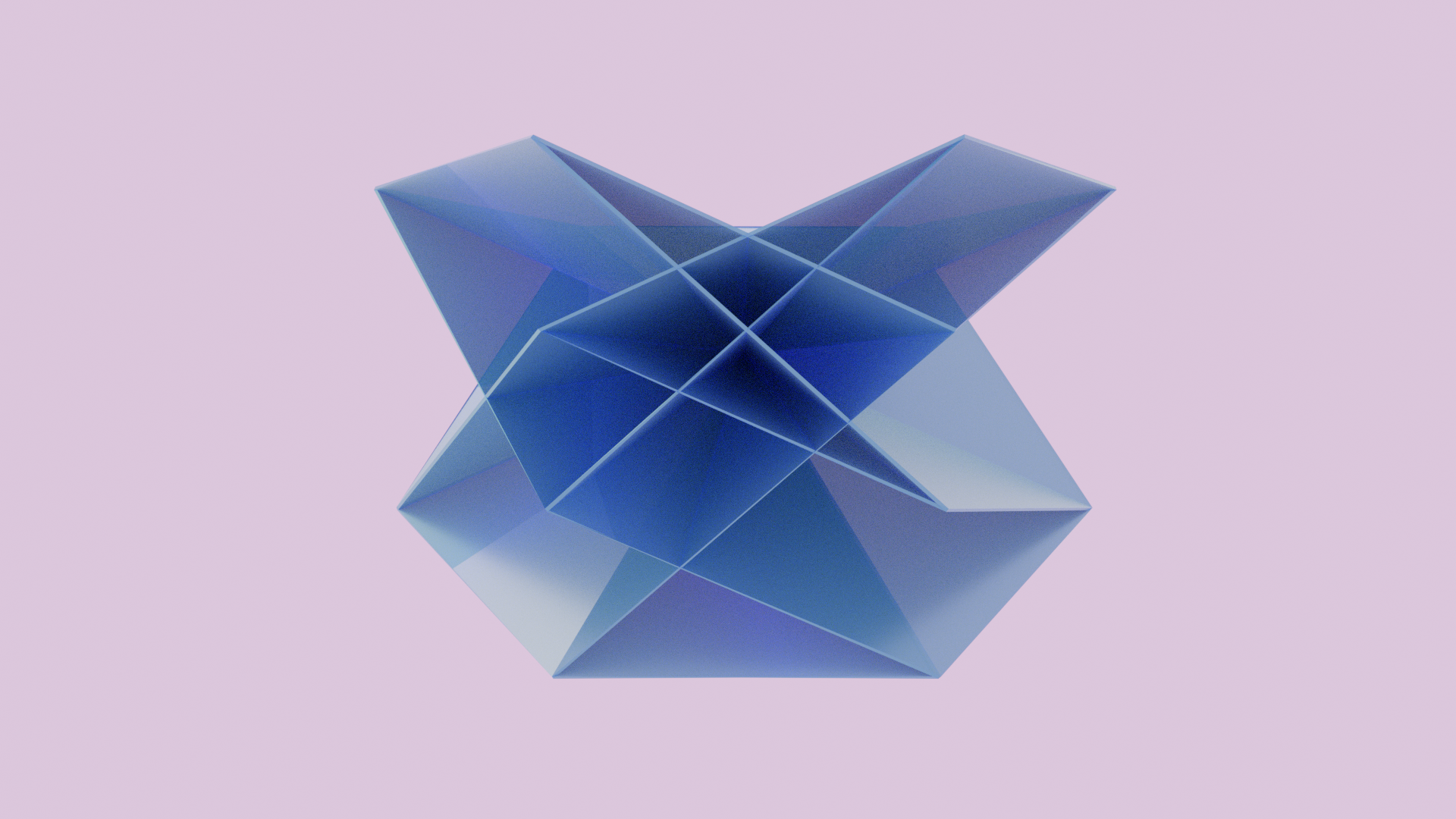}
\end{figure}

\section{Introduction}
No doubt groups are among the most versatile and omnipresent objects in mathematics, be it abstract or applied. They allow us to study mathematical constructs through their alter egos: their sets of symmetries. In important situations, this `duality' between objects and their symmetries is very precise, and offers the most important tool in our toolbox. For example in the case of linear algebra (vector spaces vs. matrix groups) or Galois theory (field extensions vs. Galois groups).

A nice method of visualising groups is that of the \emph{building} associated to a group. Given that buildings have been studied since the 1960s, introduced by Jacques Tits \citep{TitsBio2014}, it is quite surprising there have been no serious efforts yet to visualise these geometric objects in more then one dimension. Although, it should be said that their structure is very intricate, even if they are one dimensional.

Buildings are geometric objects associated to certain types of groups. They give a kind of cross section of the structure of the group, very much in the same way as Cayley or coset graphs. A building can be thought of as a higher dimensional analogue of the coset graph for a particular subgroup $B$, with the highest dimensional parts corresponding to the cosets of $B$, and lower dimensional parts corresponding to cosets of unions of conjugates of $B$.
\bigskip

\noindent The goal of the project presented here is to give a proof of concept that buildings of dimension one or two can be visualised with the help of digital methods. We consider groups with a so-called $BN$-pair in particular.

Based on existing theory, we managed to procedurally generate the structure of certain finite and infinite types of buildings using the programming language C\#, requiring minimal input variables based on the underlying group structure. Using a popular 3D-engine called Unity we visualised the buildings in a way that gives a lot of information about the underlying group.

In section 2 groups with a $BN$-pair and buildings are introduced. In section 3 we discuss the difficulties that arise in working with these groups. In section 4 we loosely explain our methods. More details can be found in the appendix.

\begin{figure}[p]
    \subfloat[The spherical building associated to the group $\text{GL}_4( \Z/2\Z)$ from the front.]{\includegraphics[width = 0.45\textwidth]{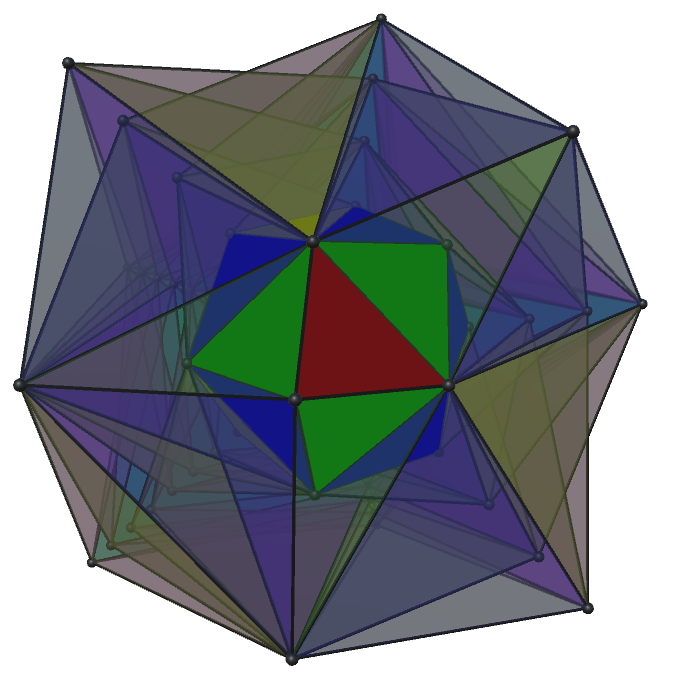}
    \label{fig:SphA3a}}\qquad
    \subfloat[The spherical building associated to the group $\text{GL}_4(\Z/2\Z)$ from the side.]{\includegraphics[width = 0.45\textwidth]{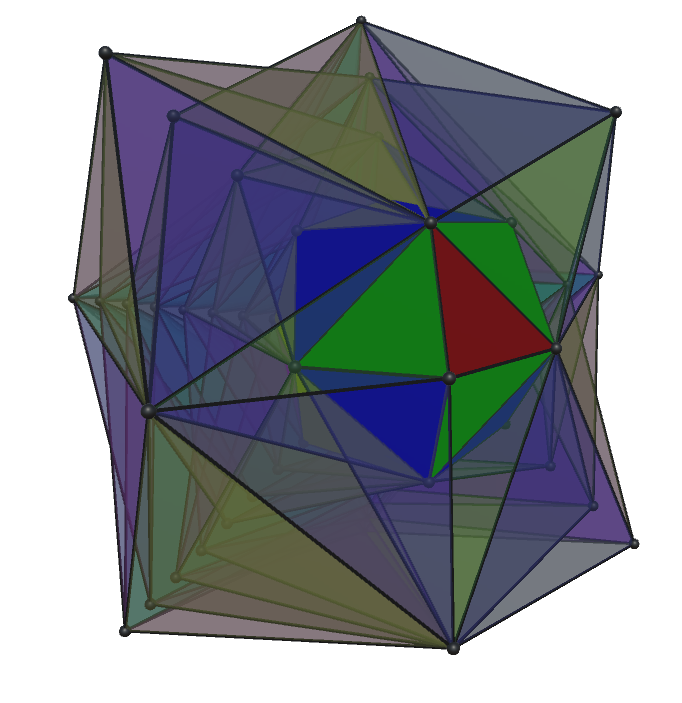}
    \label{fig:SphA3b}}
    \caption{}
    \label{fig:SphA3}
\end{figure}

\begin{figure}[p]
    \subfloat[The spherical building associated to the group $\text{GL}_3( \Z/2\Z)$.]{\includegraphics[width = 0.35\textwidth]{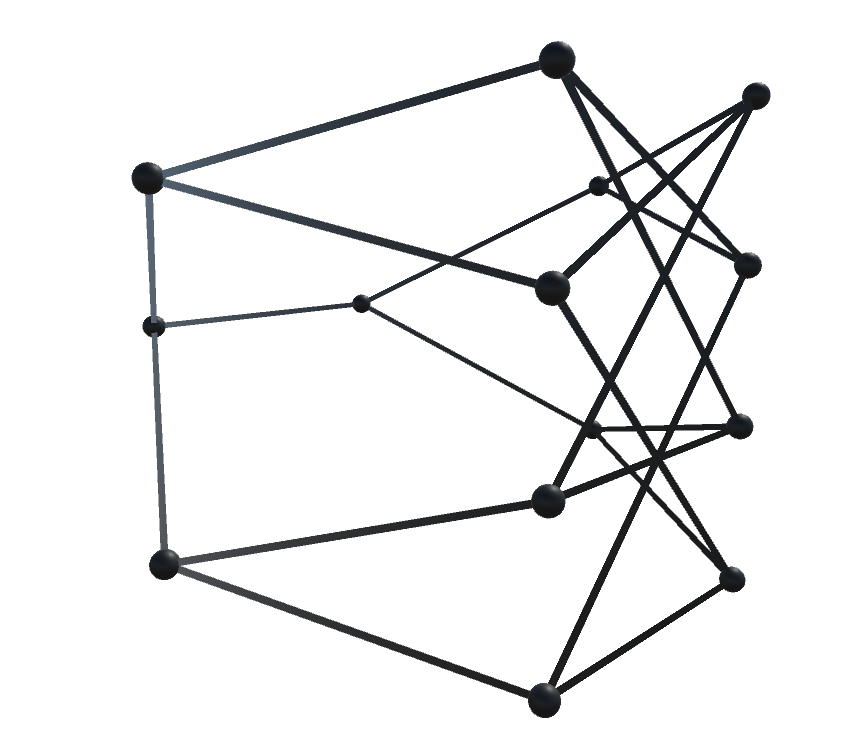}
        \label{fig:SphA2}}\qquad
    \subfloat[The affine building $\mc B(SL_3 (\Q),\nu_2)$ up to distance 5 from the chamber in the centre.]{\includegraphics[width = 0.55\textwidth]{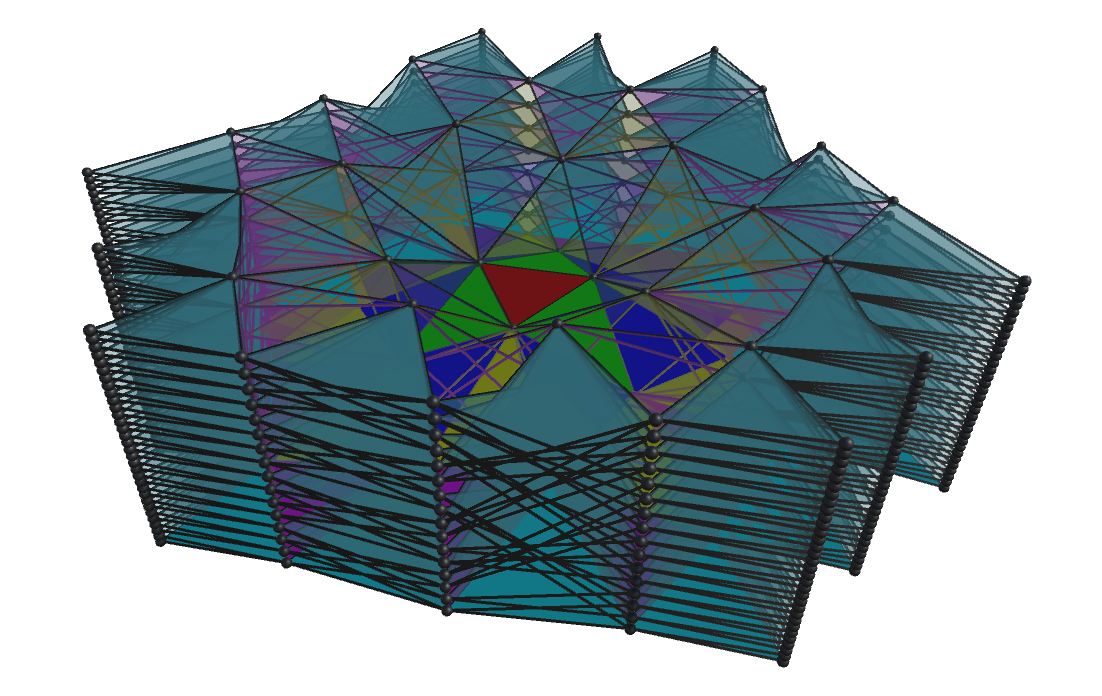}
    \label{fig:AffA2}}
    \caption{}
    \label{fig:A2}
\end{figure}

The final product is an interactive web-application, The Buildings Gallery, available at \citep{Bekker2020}. The images in Figures \ref{fig:SphA3}  and \ref{fig:A2} give a preview. We discuss our results in section 5.
\bigskip

\noindent Although the method is novel, it is entirely based on pre-existing theory. A fairly complete source for section 2 on the general theory of buildings is \citep{AbramBrown}. For specifics on $BN$-pairs we refer to \citep[\S6.2]{AbramBrown} and the original work by Bruhat and Tits \citep[\S1.2]{BrTi1}, although we keep our discussion superficial and do not point to specific results in these works.


\section{Background}

Buildings in the mathematical sense can arise in several ways, see \citep{AbramBrown,Garrett1997,Tits1997}.
We follow a group-theoretic approach, due to Tits and later worked out in collaboration with Bruhat
\citep{BrTi1}. It relies on the notion of a group with a $BN$-pair. That is a group with subgroups $B$ and
$N$ which satisfy a certain list of axioms. The precise axioms are not particularly insightful and we
do not explicitly use them in this paper, so we just refer to \cite[\S 1.2]{BrTi1}.

\begin{ex}\label{ex:BG1}
Let $F$ be a field and consider the group $G = GL_3 (F)$ of invertible $3 \times 3$-matrices with
coefficients in $F$. The subgroup
\[
B = \Big\{ \begin{pmatrix}
a & b & c \\
0 & d & e \\
0 & 0 & f \end{pmatrix} : a,b,c,d,e,f \in F, adf \neq 0 \Big\}
\]
is known as a Borel subgroup of $G$. Another important subgroup is the torus
\[
T =  \Big\{ \begin{pmatrix}
a & 0 & 0 \\
0 & d & 0 \\
0 & 0 & f \end{pmatrix} : a,d,f \in F, adf \neq 0 \Big\} .
\]
As $N$ we take the normaliser of $T$ in $G$:
\[
N = N_G (T) = \{ n \in G : n T n^{-1} = T \} .
\] 
It consists of the matrices that have precisely one nonzero entry in every column and in every row.
These $B$ and $N$ form a $BN$-pair in $G$, and one can recover $T$ as $B \cap N$. Examples of this
kind have motivated the terminology $BN$-pair. 
\end{ex}

\noindent From now on $G$ will be a group with a $BN$-pair. In first approximation, the building of $(G,B,N)$
is the set $G / B$ with the $G$-action by left multiplication:
\[
g \cdot g' B = g g' B \qquad g,g' \in G .
\]
The axioms imply that the $G$-stabiliser of $gB$ is just $g B g^{-1}$.

\begin{figure}[hbt]
    \centering
    \includegraphics[width = 0.7\linewidth]{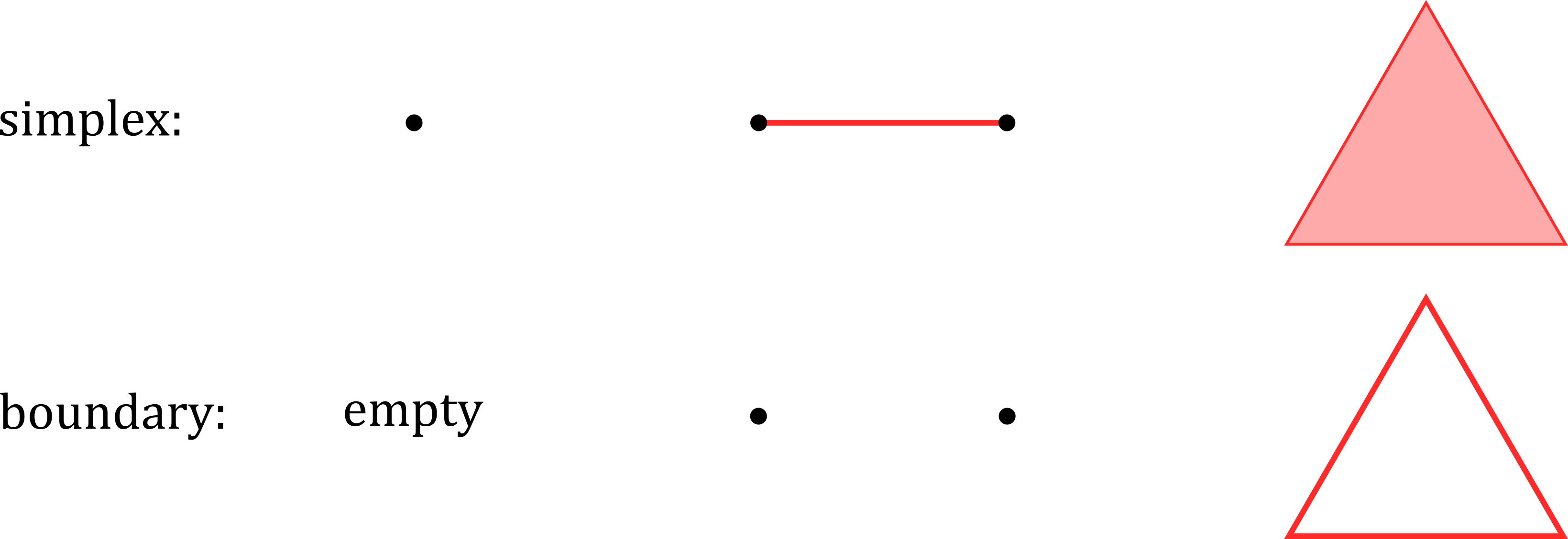}
    \caption{Some simplices and their boundaries.}
    \label{fig:boundaries}
\end{figure}

More accurately, the associated building $\mc B (G,B,N)$ is a simplicial complex: a space made from 
pieces called simplices, which are glued together in a specified way. A single simplex is a point, a 
bounded line segment, a solid triangle or a higher dimensional analogue of those. Two simplices
may only be glued along parts of their boundaries. The boundary of a bounded line segment consists
of its two endpoints, so two 0-dimensional simplices. The boundary of a solid triangle is of course
a triangle, and it can be regarded as three 1-dimensional simplices, joined along their endpoints. This is illustrated in Figure \ref{fig:boundaries}.

In the building $\mc B (G,B,N)$, a simplex of maximal dimension is called a chamber, and the set
of chambers is by definition $G / B$. All chambers have the same dimension, which is also the 
dimension of $\mc B (G,B,N)$. Every other simplex in this building can be obtained as an intersection
of a few chambers. 

\begin{figure}[hbt]
    \centering
    \includegraphics[width = 0.9\linewidth]{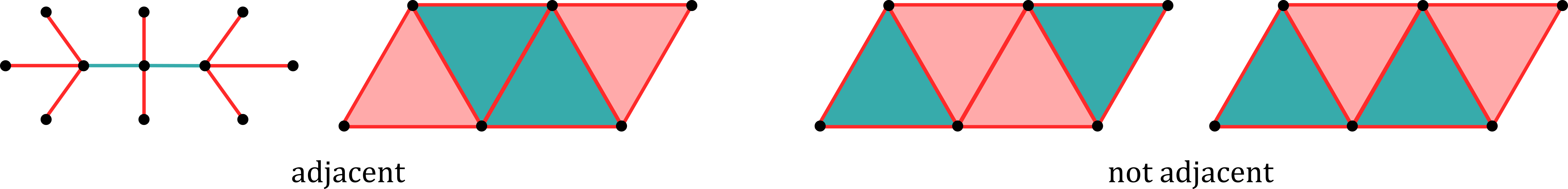}
    \caption{On the left the blue chambers are adjacent. On the right they are not.}
    \label{fig:adjacent}
\end{figure}

Two different chambers are called adjacent if their boundaries intersect in a simplex whose dimension
is one lower than that of a chamber, see Figure \ref{fig:adjacent}. In the real-estate terminology, such a codimension one simplex
in a chamber is called a wall.

By design, the building of $(G,B,N)$ is determined by the set of chambers $G / B$ and the adjacency relation on that. Adjacency can be expressed quite easily in terms of the group $G$: two different chambers $g B$, $g' B$ are adjacent if and only if the intersection of their stabiliser subgroups $g B g^{-1}$ and $g' B g^{'-1}$ is as large as possible. Equivalently, $g B$ and $g'B$ are adjacent if the group generated by $g B g^{-1} \cup g' B g^{'-1}$ is as small as possible, which means that there does not exist any subgroup $H$ of $G$ such that
\[
g B g^{-1} \subsetneq H \subsetneq \langle g B g^{-1} \cup g' B g^{'-1} \rangle .
\]
Similarly one can express group-theoretically when three (or more) chambers intersect in a simplex
of lower dimension. In this way one can construct the building $\mc B (G,B,N)$ using only group theory.

Morally speaking, this building is a geometric incarnation of the group $G$. Not only does $G$ act on
$\mc B (G,B,N)$, this group action allows one to recover and explain many interesting properties of
$G$. We note that in general a group may admit several $BN$-pairs, and hence
several buildings.

\begin{ex} \label{ex:BG2}
We consider the group 
\[
G = SL_2 (\Q) = \Big\{ \begin{pmatrix} a & b \\ c & d \end{pmatrix} : a,b,c,d \in \Q, ad - bc = 1 \Big\} .
\]
Example \ref{ex:BG1} can be adapted to this $G$, and yields a $BN$-pair with $B$ the Borel subgroup
of upper triangular matrices. But there exist other $BN$-pairs on $SL_2 (\Q)$, which are really
different. We take (essentially) the same $T$ and $N$ as in Example \ref{ex:BG1}, so 
\begin{align*}
& T = \Big\{ \begin{pmatrix} a & 0 \\ 0 & d \end{pmatrix} : a,d \in \Q, ad = 1 \Big\} , \\
& N = N_G (T) = T \cup \Big\{ \begin{pmatrix} 0 & b \\ c & 0 \end{pmatrix} : b,c \in \Q, -bc = 1 \Big\} .
\end{align*} 
To specify the new subgroup $B$, we need additional data. For a prime number $p$, the $p$-adic valuation
$\nu_p : \Q \to \Z \cup \{\infty\}$ is defined by
\[
\nu_p (q) = \left\{
\begin{array}{ll}
n & \text{if } q = p^n a / b \text{ with } n,a,b \in \Z \text{ and } a,b \text{ not divisible by } p\\
\infty & \text{if } q = 0 
\end{array} \right. .
\]
Using $\nu_p$ we define
\[
B_p = \Big\{ \begin{pmatrix} a & b \\ c & d \end{pmatrix} \in SL_2 (\Q) :
\nu_p (a), \nu_p (b), \nu_p (d) \geq 0, \nu_p (c) > 0 \Big\} .
\]
In words: the entries $a,b,d \in \Q$ do not involve powers of $p^{-1}$, while a reduced expression
for $c \in \Q$ involves at least one factor $p$.

As announced, $(B_p,N)$ forms a $BN$-pair in $SL_2 (\Q)$. We denote the associated building by
$\mc B (SL_2 (\Q), \nu_p)$. It is completely different from the buildings that arise from $BN$-pairs
like in Example \ref{ex:BG1}. A small part of $\mc B(SL_2(\Q),\nu_2)$ is drawn below, the dotted lines indicating that it continues outside the picture. For readers who know about $p$-adic numbers, we mention that
$\mc B (SL_2 (\Q), \nu_p)$ can be identified with a building associated to $SL_2 (\Q_p)$.

\begin{figure}[hbt]
    \centering
    \includegraphics[width = 0.6\linewidth]{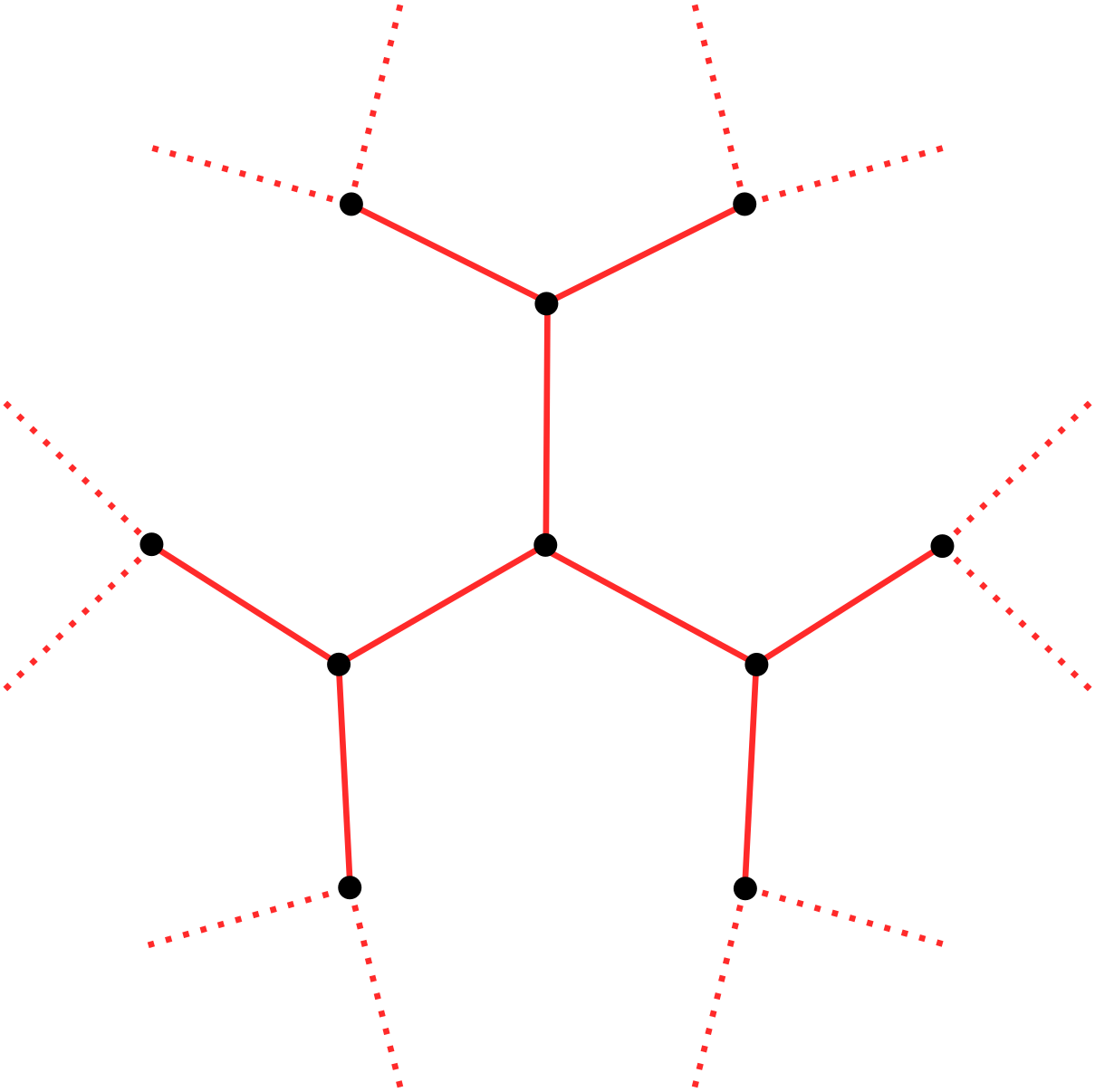}
    \caption{The building $\mc B(SL_2(\Q),\nu_2)$ is an infinite tree.}
    \label{fig:AffA1}
\end{figure}
\end{ex}
 
\noindent The buildings we encountered so far are already pretty complicated, and these are only the simplest in larger classes of examples. Although such a building always has a large symmetry group (by construction it contains $G$, but there may be many more symmetries), it is not all clear in advance how a space 
obtained by gluing chambers as above looks like. Understanding the building $\mc B (G,B,N)$ involves 
a parametrisation of the set of chambers $G / B$ by elements of $G$, so it is intimately related to 
understanding the group $G$ better. While such a parametrisation can usually be obtained in various 
ways, it tends to be really difficult to express the adjacency relations in those terms.
\bigskip

\begin{figure}[hbt]
    \centering
    \includegraphics[width = 5cm]{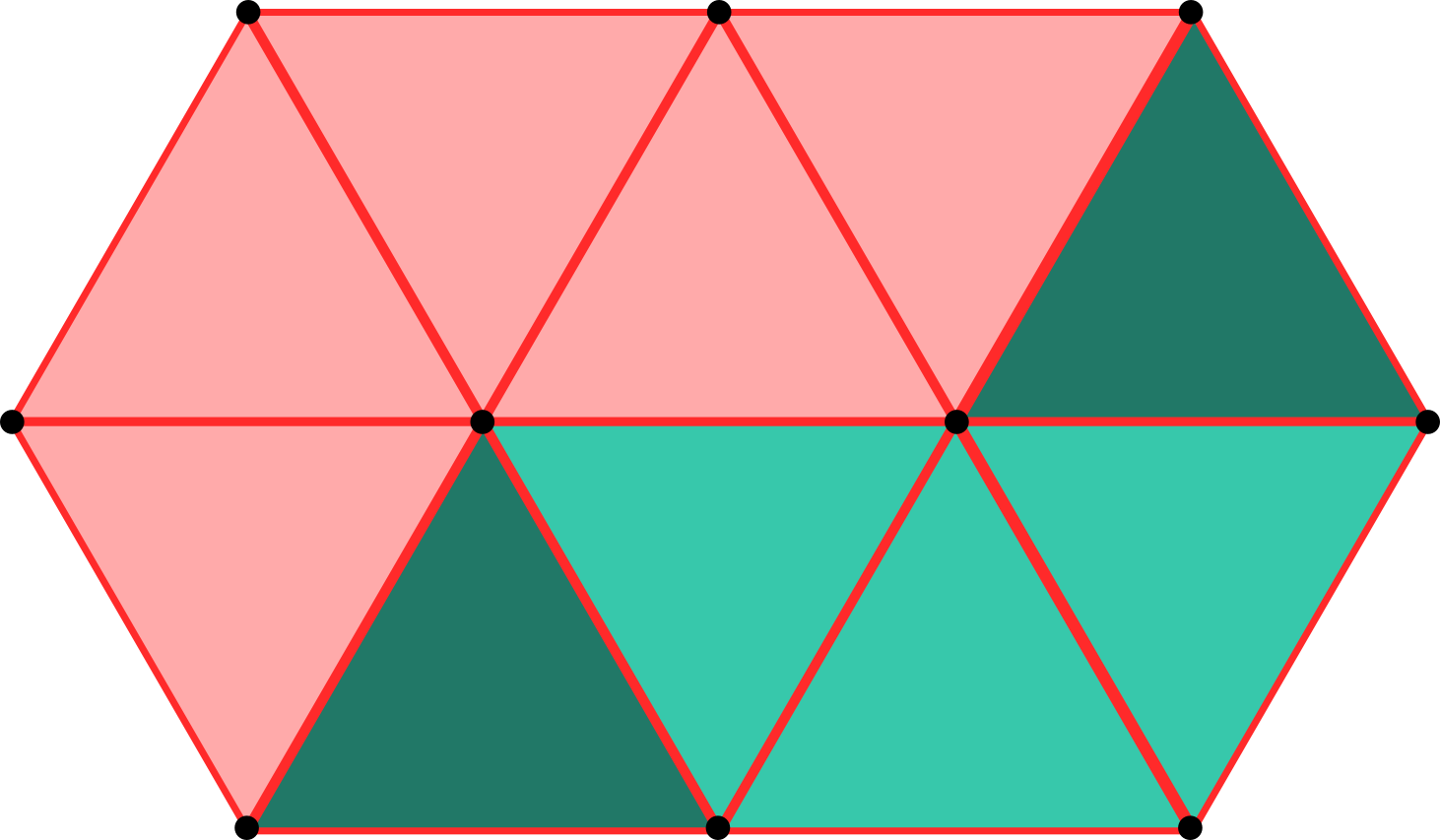}
    \caption{The distance between the dark green chambers is four.}
    \label{fig:distance}
\end{figure}

\noindent In every building, there is a useful distance function on the set of chambers. Namely, for any two
chambers $gB, g'B$ there exists a sequence of chambers
\begin{equation}\label{eq:sequence}
g B = g_0 B, g_1 B, \ldots, g_r B = g' B
\end{equation}
such that $g_i B$ and $g_{i+1} B$ are adjacent for all $i$. The length of the sequence \eqref{eq:sequence}
is $r$, and the distance between $gB$ and $g'B$ is the minimal length of a sequence with these endpoints.
For instance, in Figure \ref{fig:distance} the distance between the dark green chambers is four.
\bigskip

\noindent It may happen there is a maximum distance between any two chambers in a building. Then the building
and the underlying $BN$-pair are called spherical. The building from Example \ref{ex:BG1} is spherical,
with maximal distance 3 between chambers. This building, with $F = \Z/2\Z$, is the one in Figure \ref{fig:A2}\subref{fig:SphA2}. For a group $G$, there is usually only one spherical building
on which it acts naturally, and that will be denoted $\mc B_{sph}(G)$. Typically it comes from a 
$BN$-pair in $G$ which is constructed with a Borel subgroup like in Example \ref{ex:BG1}.

On the other hand, the building $\mc B (SL_2 (\Q), \nu_p)$ from Example \ref{ex:BG2} is not spherical. It belongs to the class of affine buildings, which means that among the non-spherical buildings it behaves as well as possible. All the buildings in this paper are either spherical or affine. An affine
building has infinitely many chambers, so it cannot be drawn entirely. Additionally, we only consider buildings that are \emph{locally finite}, meaning each chamber only has a finite number of neighbours. This allows us to focus on a part at a finite
distance from the chamber $B \in G / B$, which is guaranteed to contain finitely many chambers under this condition.

\section{Technical difficulties}
The first goal of this project was to algorithmically produce the buildings of Example \ref{ex:BG1} of the previous section with $F$ a finite field. Since the group $GL_3(F)$ is finite if $F$ is finite, we expected this to be an easy test case. The ultimate goal was to generalise the methods for this easy case, and use them to visualise the building of Example \ref{ex:BG2} of the previous section, but for $SL_3( \Q)$ instead of $SL_2( \Q )$. The building $\mc B_{sph} (GL_2(F))$ is zero-dimensional, i.e. it is just a set of points corresponding $G/B$, and the building $\mc B (SL_2(\Q),\nu_p)$ is a tree, a simple kind of graph. So the buildings $\mc \mc B_{sph} (GL_3(F))$ and $\mc B (SL_3(\Q),\nu_p )$ in this paper are the easiest more interesting examples. The building $\mc B_{sph}(GL_3(\Z/2\Z))$ is drawn in \ref{fig:A2}\subref{fig:SphA2}, and a part of $\mc B (SL_3(\Q),\nu_2)$ is drawn in Figure \ref{fig:A2}\subref{fig:AffA2}.

\begin{figure}[hbt]
    \centering
    \includegraphics[width = 0.7\linewidth]{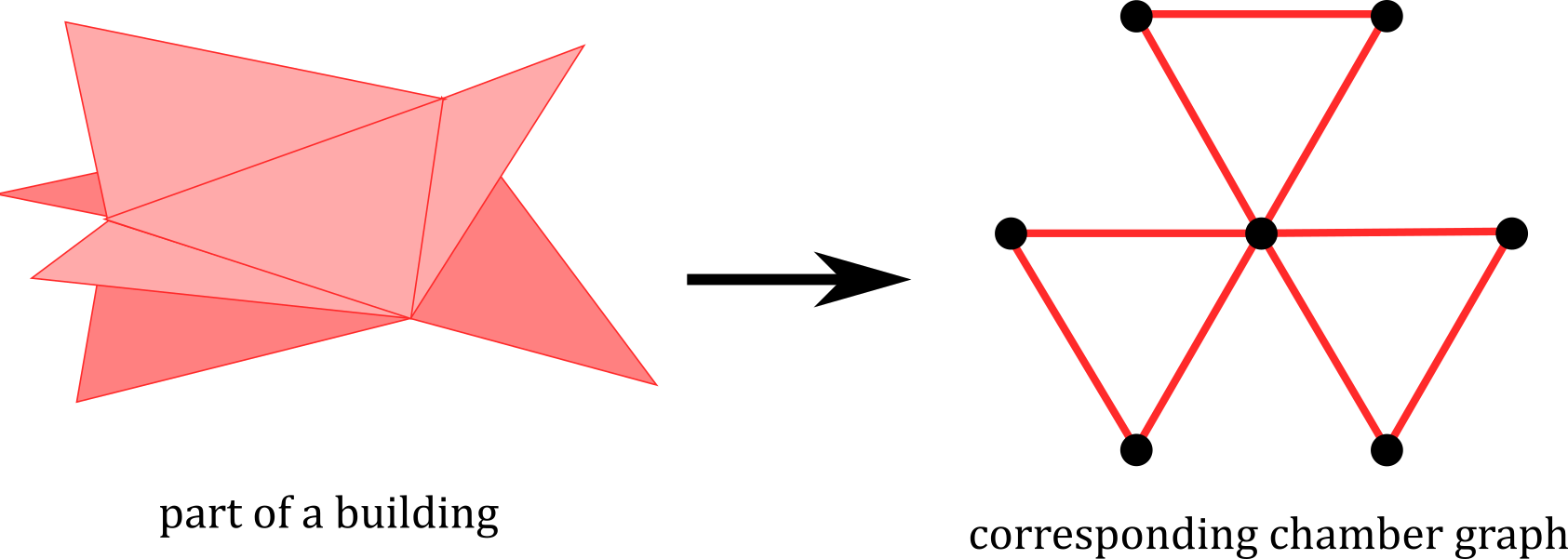}
    \caption{A part of a building with its corresponding chamber graph.}
    \label{fig:chambergraph}
\end{figure}

Recall that the chambers of the building are elements of $G / B$. To understand any building $\mc B (G,B,N)$, we have to understand $G / B$. To this end, we need to find suitable representatives of the cosets of $B$. We will refer to these representatives as the \emph{labels} of the chambers. Once we have a label for each chamber, we can use them to determine which chambers are adjacent. The set of chambers together with the set of unordered pairs of adjacent chambers is called the \emph{chamber graph}. Its vertices are the chambers, and two chambers share an edge if and only if they are adjacent. Figure \ref{fig:chambergraph} shows an example.

Obtaining a set of labels and the chamber graph is the essential problem. Once we have those, the actual visualisation does not require any new inventions. To explain why novel methods are necessary for the generation of labels and the chamber graph, we explain our early attempts first.
\bigskip

\noindent We started with our test case $\mc B_{sph} (GL_3(\Z/2\Z))$. In this paragraph we abbreviate $G = GL_3(\Z/2\Z)$. The first attempt was a naive, brute force approach: since the group $G$ is finite, it is possible to generate all of its elements, and define the subset $B$ and $N$. Finding all cosets of $B$ is then a matter of multiplying all elements of $B$ by some element $g\in G$, and checking if the resulting set contains elements of a previously generated coset. If not, it is a new coset, which we assign the label $g$.

The chamber graph can be generated by brute force as well. This is where the group $N$ comes in, which we have not used yet. Understanding why this method works requires some more advanced group theory, we refer to the appendix for a more detailed description. What is important, is that $N$ is generated by $N\cap B$ and a set $S = \{s_0, s_1\}$ of special elements we call \emph{simple reflections}. In general, $S$ can be any size. Recall that the group $N$ contains only matrices with exactly one non-zero entry in each row and column. The simple reflections for our group $G$ are
\[s_0 = \begin{pmatrix} 0   &1  &0\\
                        1   &0  &0\\
                        0   &0  &1  \end{pmatrix},\ 
  s_1 = \begin{pmatrix} 1   &0  &0\\
                        0   &0  &1\\
                        0   &1  &0  \end{pmatrix}.\]
They are called reflections because in the geometric realisation of the building they act like reflections. It turns out that two chambers $gB$ and $hB$ are adjacent if and only if there exists an $s_i\in S$ and $b,b'\in B$ such that $h =gbs_ib'$. The sets involved are again all finite, and small enough to make checking this for every combination of $s_i,b'$ and $b$ feasible.
\bigskip

\noindent Of course, these methods are very inefficient. They only work well because the group $\text{GL}_3\left(\mathbb{Z}/2\mathbb{Z}\right)$ is very small. However, both the group $SL_3 ( \Q)$ and the building $\mc B (SL_3 ( \Q),\nu_p)$ of Example \ref{ex:BG2} are infinite. There are in fact three problems caused by infinite sets:
\begin{itemize}
    \item First of all, the set of chambers $G / B$ might in general be infinite.
    \item Secondly, the subgroup of $N$ generated by the set of simple reflections $S$ might be infinite, even if $S$ contains only a couple of elements.
    \item Finally, the subgroup $B$ might be infinite.
\end{itemize}
The first two problems are closely related. For our examples, they have to do with the fact that the distance between two chambers can be arbitrarily large, which is the case for $\mc B(SL_3 (\Q),\nu_p)$. In any of these cases there is no hope of a general method by brute force.

At this point, a seasoned programmer would probably start looking for a work-around. There might be a way to achieve the same result, without having to do all the work. For example, is it really necessary to calculate the entire chamber graph? A useful property of buildings is that their structure is very uniform, in the sense that every `local part' of the building looks roughly the same. Could we not `copy-paste' certain finite parts?

Maybe this could work. For example, these `local parts' of the building $\mc B (SL_3(\Q),\nu_p)$ look like the spherical building $\mc B_{sph} (GL_3 (\Z/2\Z))$. So we could streamline our brute-force approach a little, and apply it to every such local part. However, it would require keeping track of how all of these local parts glue together, and that is not a trivial task either. Besides, that would not solve the issue of the infinitely many cosets of $B$. The use of some proper group theory to solve these issues is inevitable.

\section{The method}
\begin{wrapfigure}[22]{r}{0.4\textwidth}
    \centering
    \includegraphics[width = 0.38\textwidth]{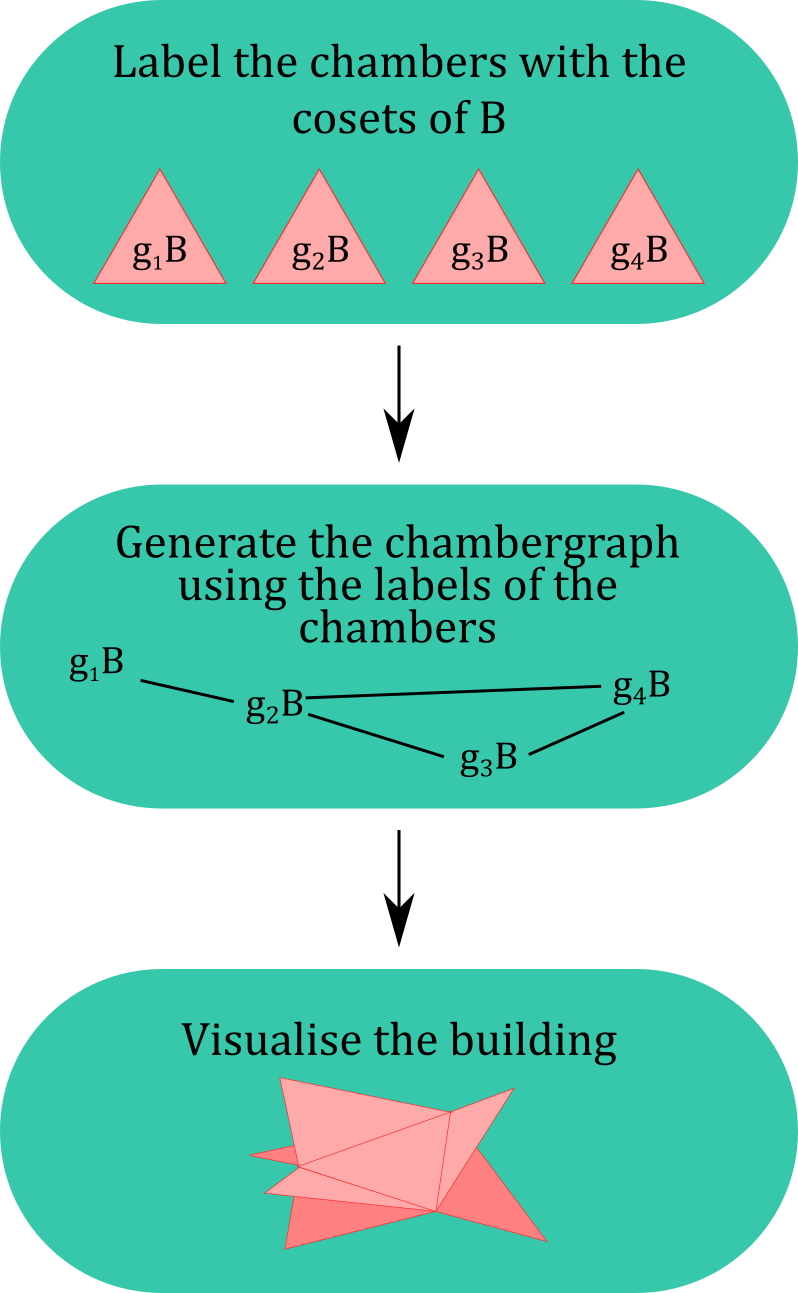}
    \label{fig:flowchart}
\end{wrapfigure}
Our solution consists of three steps: generate labels for the chambers in terms of representatives of $G/B$, generate the chamber graph and visualise the building -- see the flowchart to the right.
The essence of the problem is formed by the first two parts, for which there were no procedural methods yet. In our approach to the first step, we first solve the issue of possibly infinitely many chambers, which immediately solves the issue of the group generated by $S$ being possibly infinite as well. In dealing with these difficulties, it also becomes relevant that we handle $B$ in a suitable way. This turns out to involve some rather advanced group theory, which we discuss in more depth in the appendix. In the second step, the group $B$ pops up again, and we need to deal with it in yet another way, however this time the solution is more straightforward.

In all that follows, $G$ is a group with a $BN$-pair. We only consider those for which $S$ is finite. In all our examples, $S$ will only have two or three elements. We assume that the building of $G$ is either spherical or affine, and locally finite. If it is spherical, it contains only finitely many chambers, and we assume the group $G$ is similar to Example \ref{ex:BG1}. If it is affine, it has infinitely many chambers, and the $BN$-pair is similar to that in Example \ref{ex:BG2}.

\subsection{The labels}
We first have to bound the number of chambers we consider. This has a simple fix: we only consider a finite part of the building, say up to a distance $d$ of the chamber $B$ in all directions. Let $gB$ and $hB$ be two chambers. Just like we have a criterion for adjacency in terms of the simple reflections, we have a criterion for $gB$ and $hB$ being any distance apart. For example, if $S = \{s_0,\ldots,s_n\}$ is the set of simple reflections, $gB$ and $hB$ are a distance 2 apart if and only if there is some $n=s_is_j$ with $i\not=j$ such that
\begin{equation}\label{eq:M1}
\text{there exist } b,b' \in B \text{ such that } h = g b n b'  .
\end{equation}
This $n$ denotes a path through consecutively adjacent chambers:
\[gB = h_0B,\ h_1B,\ h_2B = hB,\]
where $h_1 = h_0b_0s_ib_0'$ and $h_2 = h = h_1b_1s_jb_1'$, with $b_i,b_i'\in B$.

\begin{figure}[hbt]
    \centering
    \includegraphics[width = 0.4\linewidth]{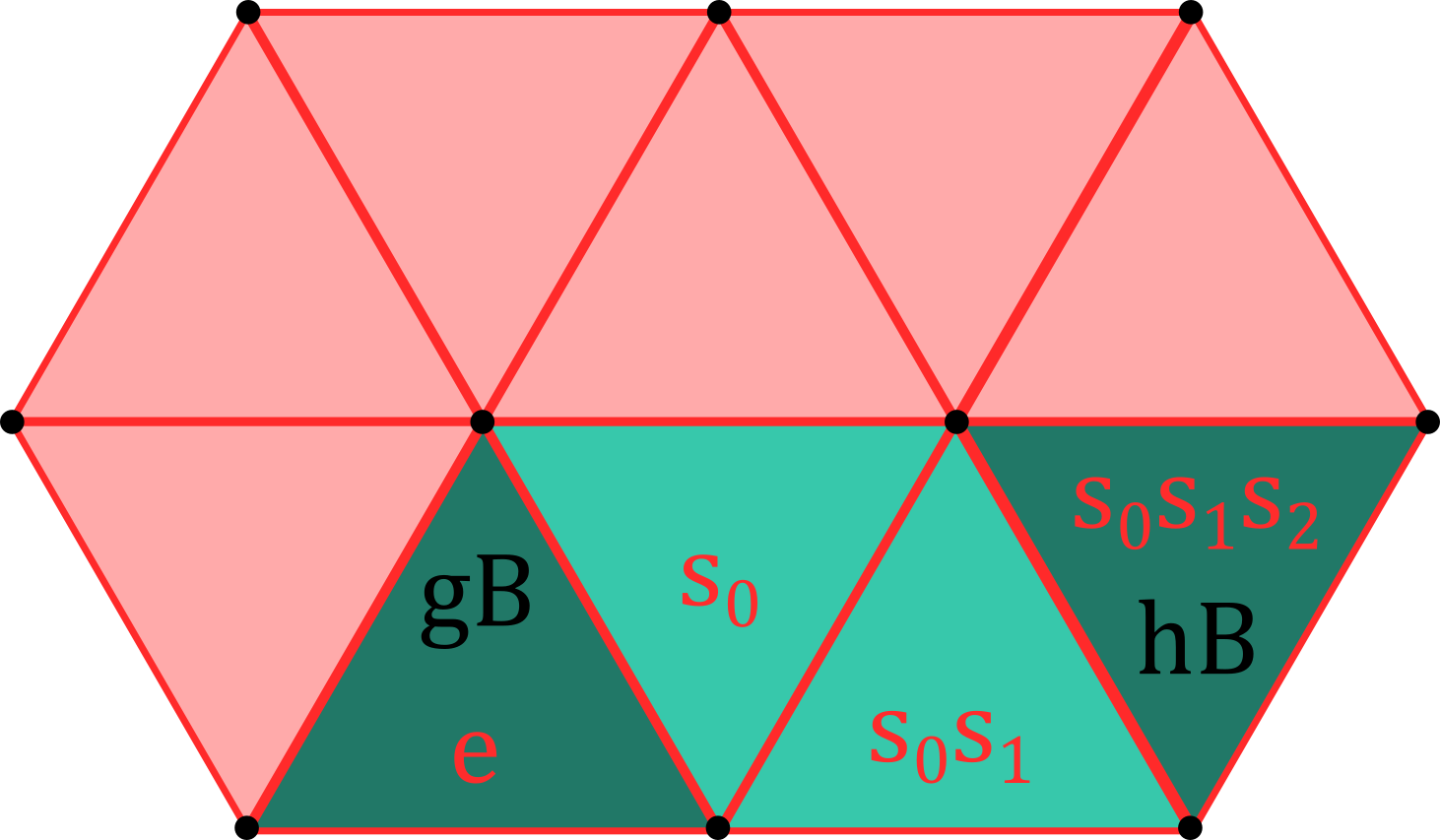}
    \caption{The shortest path between $gB$ and $hB$ in terms of some element $n\in N$. The element $n$ of equation \eqref{eq:M1} for a chamber on this path is written in red. The unit element is denoted by $e$.}
    \label{fig:wdistance}
\end{figure}

In general, two chambers $gB$ and $hB$ are distance $r$ apart if and only if \eqref{eq:M1} holds with $n$ equal to a product $s_{i_1}\cdots s_{i_r}$, and $r$ is the smallest number of $s_i$ needed to write $n$ as such a product. This is demonstrated in Figure \ref{fig:wdistance} below for a building with three simple reflection $s_0,\ s_1$ and $s_2$. We call $r$ the \emph{length} $l(n)$ of $n$. If we choose $g=e$, where $e$ is the unit element, any other coset associated to a chamber we consider then takes the form
\[bnB,\text{ with }b\in B\text{ and }n\in N:\ l(n)\leq d.\]
Since $S$ is finite, there are only finitely many possibilities for $n$, which all give distinct cosets.

This still leaves possibly infinitely many options for $b$. However, for the buildings that we consider, for any fixed choice of $n$ the number of cosets $bnB$ with $b\in B$ is finite. So, we should be able to write down a finite subset of $B$ parametrising the chambers for $n$. This is the most technical part of our approach. In the appendix we derive an algorithm for finding these representatives for any fixed $n$.

\subsection{The chamber graph}
Labelling the chambers with elements of $G$ enables us to decide whether two chambers $gB$ and $hB$ are adjacent. We have to check whether $h = gbs_ib'$ for some $s_i\in S$ and $b,b'\in B$. In the appendix we give the interested reader some idea of where this criterion comes from. In general, finding a decomposition of $h$ in these terms can be very hard, if not unsolvable. An alternative is to find a decomposition like this by hand, as we did before. Unfortunately, in general we would have to check this condition for infinitely different $b\in B$.

However, we found that we do not have to make many calculations; the elements of $Bs_iB$ have certain restrictions on their coordinates, which are distinct for each $s_i$. We only have to check whether the coordinates of the matrix $g^{-1}h$ meet these requirements. Since there are at most three $s_i$, we can explicitly hard code these requirements without much effort.

\begin{ex}
    Let $G = GL_3 (\Z/2\Z)$, and recall that $B$ is the set of invertible upper triangular matrices, and
    \[s_0 = \begin{pmatrix} 0 & 1 & 0\\
                            1 & 0 & 0\\
                            0 & 0 & 1\end{pmatrix},
      s_1 = \begin{pmatrix} 1 & 0 & 0\\
                            0 & 0 & 1\\
                            0 & 1 & 0\end{pmatrix}.
                            \]
    Then, the sets $Bs_0B$ and $Bs_1B$ are given by
    \[\begin{split}
        Bs_0B   &= \left\{\begin{pmatrix}   a & b & c\\
                                            1 & d & e\\
                                            0 & 0 & 1\end{pmatrix}\in GL_3(\Z/2\Z)\right\},\\
        Bs_1B   &= \left\{\begin{pmatrix}   1 & a & b\\
                                            0 & c & d\\
                                            0 & 1 & e\end{pmatrix}\in GL_3(\Z/2\Z)\right\}.
    \end{split}
    \]
    \label{ex:M1}
\end{ex}

\begin{ex}\label{ex:M2}
    Let $G = SL_3(\Q)$ and consider the $BN$-pair corresponding to $\mc B (SL_3( \Q),\nu_p)$. It is given by matrices $b$ whose coordinates in the upper triangular part have $p$-adic valuation $\geq 0$, in the lower triangular part have $p$-adic valuation $>0$, and on the diagonal have $p$-adic valuation exactly 0. In other words: $b\in B_p$ is such that the coordinate $b_{i,j}$ has
    \[\nu_p(b_{i,j}) \begin{cases}    \geq 0\text{, if }i<j\\
                                        = 0\text{, if }i=j\\
                                        > 0\text{, if }i>j\end{cases}.\]
    
    In comparison with the spherical $BN$-pair for $SL_3 (\Z/2\Z)$, there is an additional simple reflection. The simple reflections are:
    \[s_0 = \begin{pmatrix} 0 & 1 & 0\\
                            1 & 0 & 0\\
                            0 & 0 & 1\end{pmatrix},\ 
      s_1 = \begin{pmatrix} 1 & 0 & 0\\
                            0 & 0 & 1\\
                            0 & 1 & 0\end{pmatrix},\ 
      s_2 = \begin{pmatrix} 0 & 0 & -p^{-1}\\
                            0 & 1 & 0\\
                            p & 0 & 0\end{pmatrix}\]
    We can then calculate the valuations of the coordinates of an arbitrary element of $B_ps_iB_p$. They are:
    \[\begin{split}
        B_ps_0B_p   &= \left\{\begin{pmatrix}   \nu_p \geq 0 & \nu_p \geq 0 & \nu_p \geq 0\\
                                            \nu_p = 0 & \nu_p \geq 0 & \nu_p \geq 0\\
                                            \nu_p \geq 1 & \nu_p \geq 1 & \nu_p = 0\end{pmatrix}\right\},\\
        B_ps_1B_p   &= \left\{\begin{pmatrix}   \nu_p = 0 & \nu_p \geq 0 & \nu_p \geq 0\\
                                            \nu_p \geq 1 & \nu_p \geq 0 & \nu_p \geq 0\\
                                            \nu_p \geq 1 & \nu_p = 0 & \nu_p \geq 0\end{pmatrix}\right\},\ \text{and}\\
        B_ps_2B_p   &= \left\{\begin{pmatrix}   \nu_p \geq 0 & \nu_p \geq 0 & \nu_p = -1\\
                                            \nu_p \geq 1 & \nu_p = 0 & \nu_p \geq 0\\
                                            \nu_p \geq 1 & \nu_p \geq 1 & \nu_p \geq 0\end{pmatrix}\right\}.
    \end{split}\]
\end{ex}

\noindent Recording all adjacency relations between chambers provides the information we need for generating the chamber graph. In fact, we obtain more. Not only do we have a set of vertices (the chambers) and edges (pairs of adjacent chambers), but for each adjacent pair we also found the simple reflections $s_i$ by which the chambers are adjacent. We include this in the chamber graph, giving each edge a label corresponding to an element of $S$. By performing a search algorithm, for example Dijkstra's algorithm, on this labelled graph, we can not only find a numerical distance between chambers $gB$ and $hB$, but also find the element $n\in N$ that determines the path through the building that brings us from $gB$ to $hB$. This is important in the next step.

\subsection{Visualisation}
This final step needs a little introduction. We called the elements of $S$ simple reflections. For good reason: they act as reflections in the geometric sense. Remember that geometrically a chamber is a simplex of highest dimension in the building. Two chambers $gB$ and $hB$ are adjacent if and only if their intersection is a wall, and if and only if $g^{-1}h \in Bs_iB$, for some $s_i\in S$. We can associate to every wall of the chamber $gB$ a simple reflection $s_i$, and think of each adjacent chamber as being obtained by reflecting $gB$ in that wall. In general $gB$ will have multiple adjacent chambers sharing the same wall, as there might be distinct $bs_iB\not=b's_iB$ for $b,b'\in B$ and $b\not=b'$. We call the reflection $s_i$ corresponding to a wall the \emph{type} of that wall.
\bigskip

\noindent We need a couple of objects. First, for each chamber for which we have generated a label, the program needs an object that holds a few attributes. These are at least the label and a set of vertices that define the simplex associated to the chamber. Since we allow $\#S$ to be two or three, the number of walls of such a simplex is two or three, so it will be either a line segment or a triangle. We give the object that we associate to a chamber also a set of walls, again defined by the vertices they contain, and we give the walls their type $s_i$.

When two chambers share walls or vertices, we need to make sure that the associated objects do so too. To this end, we can just look at what simple reflection they are adjacent for, find the wall of the corresponding type and copy the vertices defining that wall and the wall itself from one chamber to the other.
\bigskip

\noindent The positions of the walls and chambers in the visualisation are entirely determined by the positions of their vertices. So all we have to do is find out where the vertices go. The reflections offer a straightforward way of doing this. Any chamber can be reached from $B$ by consecutively reflecting in some sequence of walls, and so we only have to hard code the position of the vertices of $B$. Then, the position of vertices of a chamber sharing a wall of type $s_i$ is obtained by reflecting the vertices of $B$ not contained in the neighbouring chamber in the line that the wall lies on. See Figure \ref{fig:reflection} below.
\begin{figure}[hbt]
    \centering
    \includegraphics[width = 0.35\linewidth]{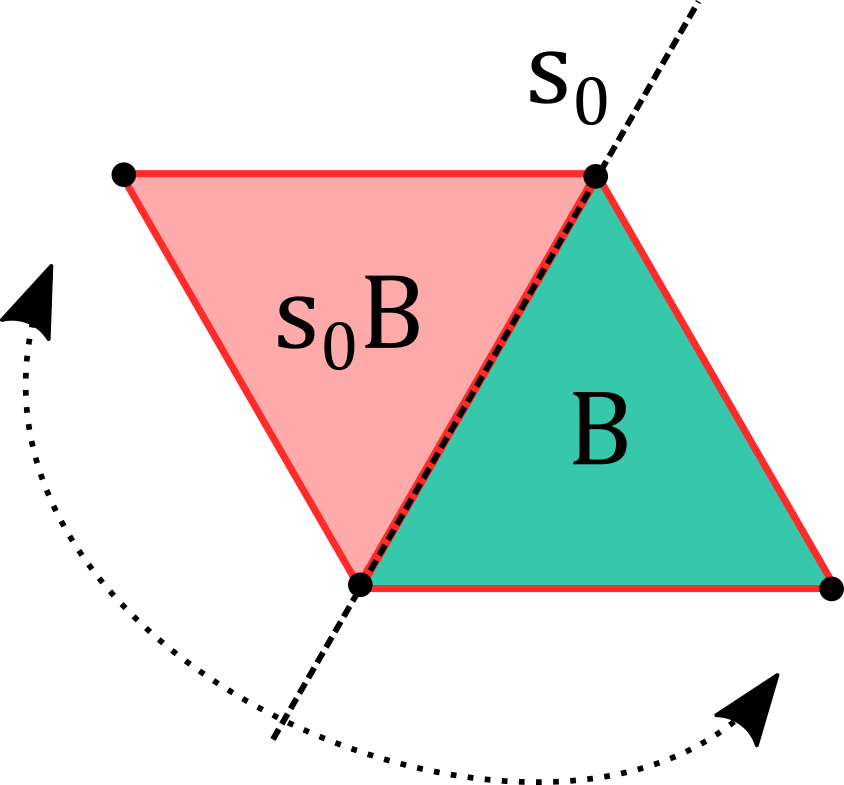}
    \caption{The positions of the vertices of chamber $s_0B$ are obtained by reflection in the wall of type $s_o$.}
    \label{fig:reflection}
\end{figure}

\begin{figure}[hbt]
    \centering
    \subfloat[A larger part of the building $\mc B(SL_3(\Q),\nu_2)$. The elements $s_i$ determine the horizontal position of the chambers.]{
    \includegraphics[width = 0.45\textwidth]{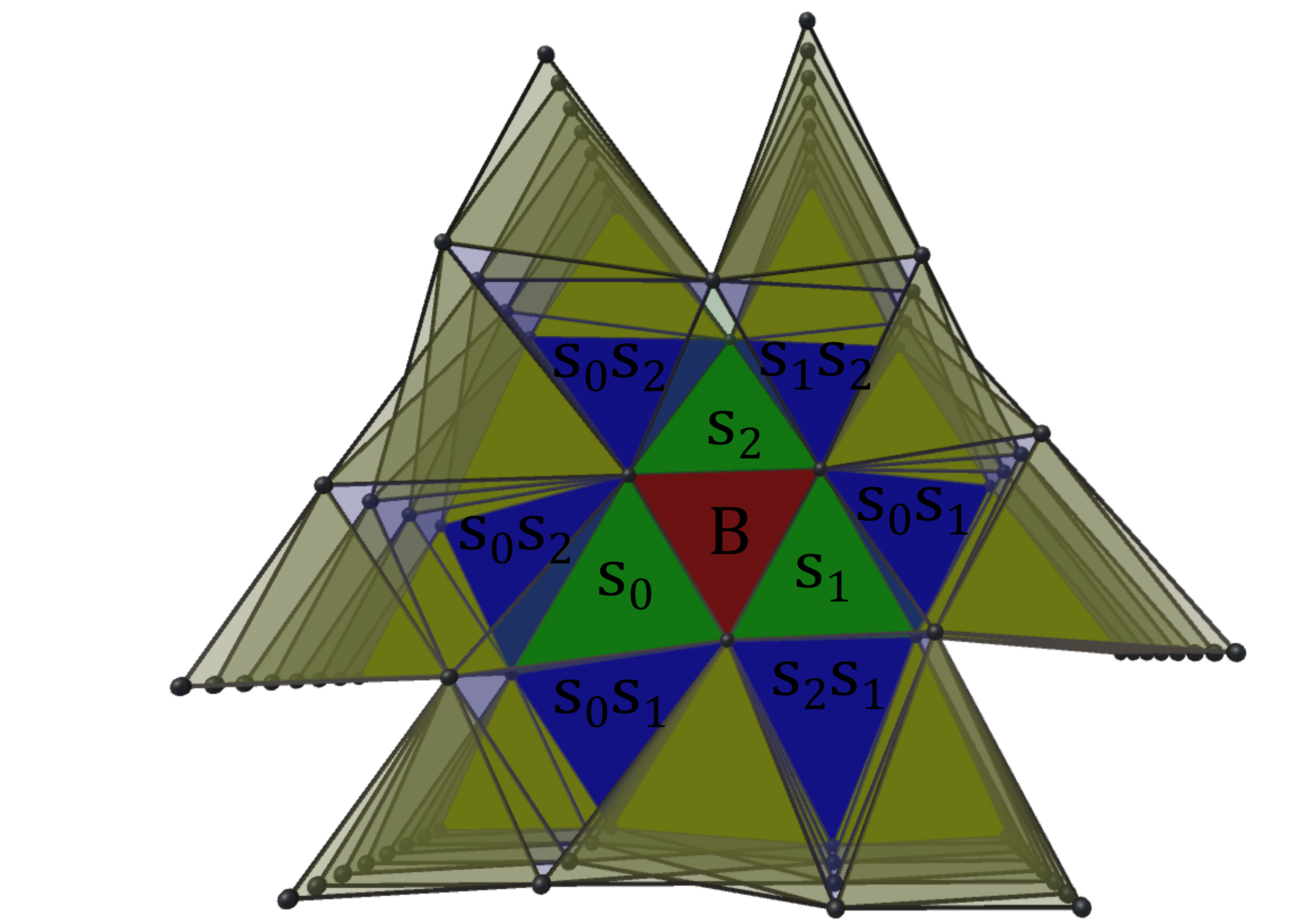}
    \label{fig:AffVisa}}\qquad
    \subfloat[A side view of the same part of the building. The height of a chamber in a `stack' of chambers related to $B$ by the same $n\in N$ is determined by the order in which the chambers are generated.]{
    \includegraphics[width = 0.45\textwidth]{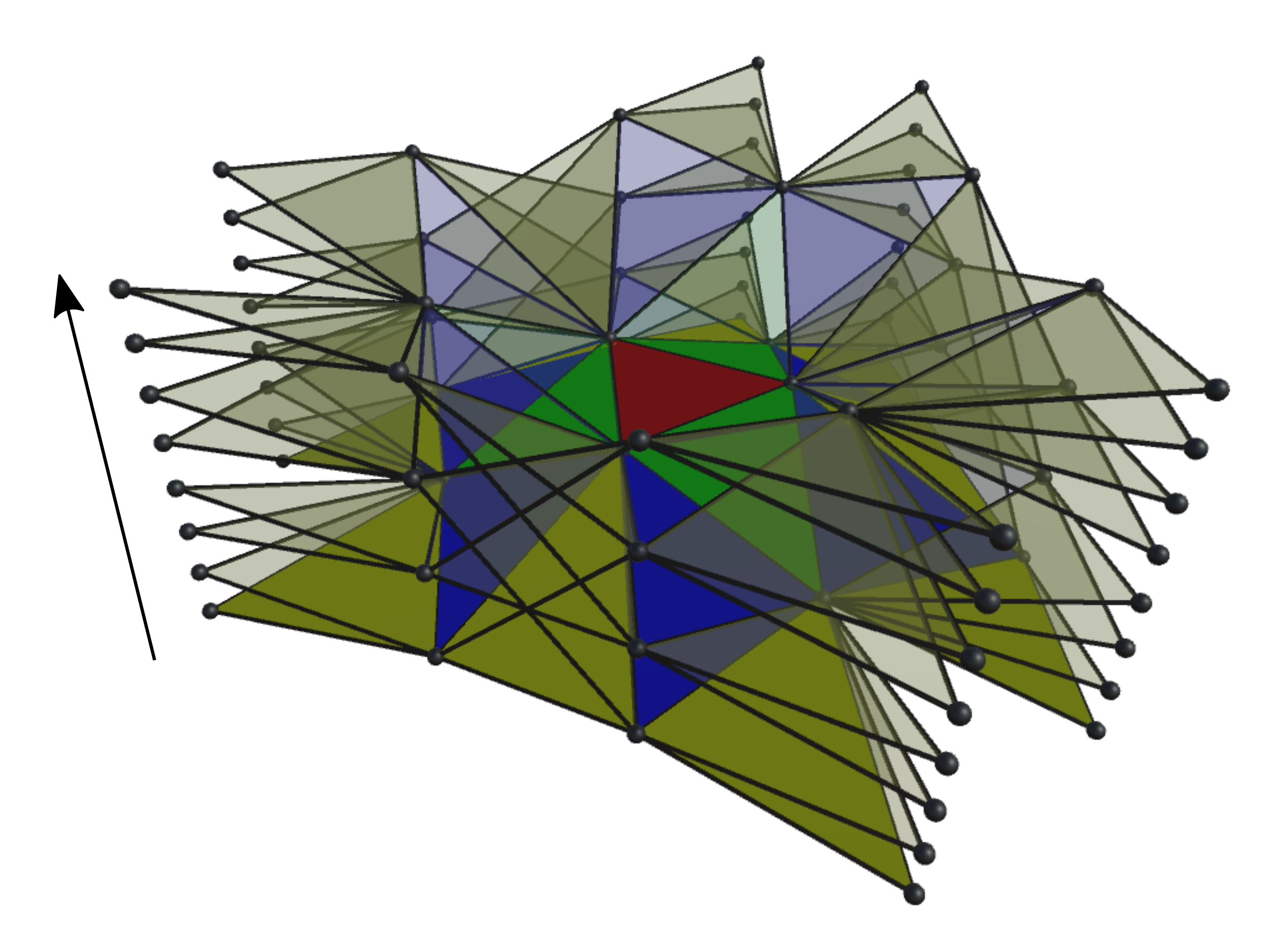}
    \label{fig:AffVisb}}
    \caption{}
    \label{fig:AffVis}
\end{figure}

In general, we can find a minimal path between two chambers $B$ and $gB$ through the chamber graph. This minimal path is represented by an element $n\in N$, which is the product of the consecutive simple reflections that takes $B$ to $gB$. If $n = s_{i_1}\cdots s_{i_r}$ for $s_{i_j}\in S$ we determine the position of the vertices of $gB$ by consecutively reflecting the vertices of $B$ in the lines containing the walls corresponding to the $s_{i_j}$. This is demonstrated in Figure \ref{fig:AffVis}\subref{fig:AffVisa}.

This is not all, since all cosets of $B$ related to it by the same element $n\in N$ now end up in the exact same location in the plane at height 0. To resolve this we separate vertices based on the order in which the chambers containing them were generated. To each chamber we assign a \emph{height}. In the case $\mc B(SL_3(\Q),\nu_p)$ we add the height of the chamber closest to $B$ containing the vertex, which is unique, to the vertical coordinate of the vertex. See Figure \ref{fig:AffVis}\subref{fig:AffVisb}.

In the building of Figure \ref{fig:SphA3} the chambers form a sphere, and so this method of separating them does not produce a nice image. Instead, we separate them radially, again based on the height we assigned each chamber. These are the only two cases we encounter when visualising spherical or affine buildings of dimension two or three.

\section{Results}
We have in total visualised three buildings using this method: our test case $\mc B_{sph} (GL_3(\Z/2\Z))$ of Figure \ref{fig:A2}\subref{fig:SphA2}, the affine building $\mc B( SL_3(\Q), \nu_2)$ of Figure \ref{fig:A2}\subref{fig:AffA2} and an additional spherical building $\mc B_{sph} (GL_4(\Z/2\Z))$, see Figure \ref{fig:SphA3}. The main difference between the latter and the first is that the latter is of dimension two, whereas the building of $GL_3(\Z/2\Z)$ is of dimension one. Our method is sufficiently general that it can handle these small differences without problem.

The building of $\mc B_{sph} (GL_3(\Z/2\Z))$ is small enough so that we can calculate its structure with purely algebraic methods by hand. This is a standard exercise, see \cite[Exercise 4.23]{AbramBrown}. Therefore, its structure has been known for a long time, and is more commonly known as the \emph{incidence graph} of the Fano plane, \cite[Fig 4.1]{AbramBrown}. The only difference between the cited figure and our visualisation in Figure \ref{fig:A2}\subref{fig:SphA2} is the choice of position of the vertices, which emphasises different aspects of the symmetry group of the building. This is a non-essential difference, and we find that our method reproduces the exact same building.

Although the other spherical building also corresponds to a finite group, it is a lot bigger, and we are unaware of any published calculation of its structure. However, it is in general verifiable whether two chambers should indeed be adjacent by a calculation on their labels, as explained in the method section. As far as we can see, the produced images are consistent with these calculations.
\bigskip

\begin{figure}[hbt]
    \centering
    \subfloat[The spherical building $\mc B_{sph}(GL_3(\Z/2\Z))$.]{\includegraphics[width = 0.4\textwidth]{SphA2Full.png}
    \label{fig:SphA2Fulla}}\qquad
    \subfloat[The spherical building $\mc B_{sph}(GL_3(\Z/2\Z))$ drawn as a part of $\mc B(SL_3(\Q),\nu_p)$. Here we include the simplices of $\mc B(SL_3(\Q),\nu_p)$ that have non-empty intersection with it.]{\includegraphics[width = 0.4\textwidth]{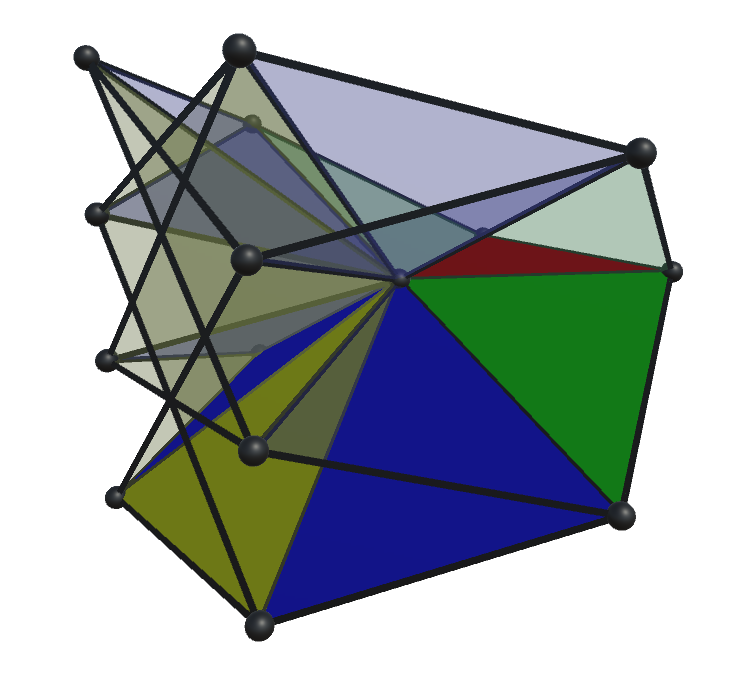}
    \label{fig:SphA2Fullb}}
    \caption{}
    \label{fig:SphA2Full}
\end{figure}

\noindent Each spherical building forms a `local' part of an affine building, a fact that we have alluded to before in this text. It is a small effort to show how the building $\mc B_{sph}(GL_3 (\Z / 2 \Z))$ would fit in the building $\mc B (SL_3(\Q), \nu_2)$, by adding one extra vertex in the centre. In Figure \ref{fig:SphA2Full} we show this building again, and the building with this vertex added. It is insightful to compare this to Figure \ref{fig:A2}\subref{fig:AffA2} by identifying the blue triangle of the spherical building with the blue triangle in the affine building.
\bigskip

\noindent The end product of this project is an online application called \emph{The Buildings Gallery}. It can be found at https://buildings.gallery\citep{Bekker2020}, and displays interactive models of the buildings. The application is under continuous development; the plan is to keep adding functionalities in the future.

\section{Outlook}
There are many things still to be explored in the topic of these visualisations. The three buildings we have visualised are only a very limited selection of what is out there. In particular, there are also buildings whose chambers are not simplices, but for example squares or cubes. The question arises whether our method is general enough to handle these cases as well.

To further illustrate certain properties, animations and other functionalities can be added to the web application. It would be nice to be able to show the possible shortest paths between chambers, or animate the action of the subgroup $B$ on the building. These are only examples of what is possible with our tool.

Finally, there are a few points on which the method can be improved and further generalised. For example, at this stage we have to hard code what the elements in the sets $BsB$ look like. It would be useful to find an algorithm to calculate these too, so that the procedure truly only depends on the groups $B$, $N$ and $T$.

That said, the produced 3D models provide a beautiful and useful visualisation of a core subject in mathematics. As such, they form an interesting topic for science communication for broad audiences, for example as maths-art. Besides this, the resulting application is a useful tool for teaching general building theory, and might inspire further research into the subject.

\section*{Acknowledgements}
We are grateful to the George Mason University (Fairfax, USA), in particular to Dr. Anton Lukyanenko and all students and staff at the Mason Experimental Geometry Lab. Thanks to them, working on this project was a delight. Dr. Lukyanenko's support and advice were essential to the success of the project.

Gratefulness also goes out to the Radboud Honours Academy, that made this project possible.

\section*{Funding}
This work was supported by the Beyond the Frontiers programme of the Radboud Honours Academy through a personal scholarship.

\section*{Disclosure statement}
The authors declare no conflict of interest.

\section{Appendix}
This appendix shows some of the inner workings of our methods, using more advanced mathematics. Readers familiar with the basics of algebraic groups and their root systems over discretely valued fields should be able to understand the method entirely, but it might be interesting for readers with less specific knowledge as well. We refer to \cite[Ch. 6]{AbramBrown} and the introduction of \citep{SchneiStuhl} for most details in this appendix. In all that follows, $G$ is a group with a $BN$-pair.

\subsection{Notes on the criterion for distance}
A central assertion in section 4 is that the distance between two chambers $gB$ and $hB$ can be measured in terms of an element of $N$. We will not rigorously derive this statement, but we will try to shine some light on it. Any group $G$ with a $BN$-pair admits a so-called \emph{Weyl group}. We first define a normal subgroup of $N$ by $T = B\cap N$. The axioms of a $BN$-pair entail that the quotient
\[W = N/T\]
is a Weyl group. An important point is that it is generated by the image of the set of simple reflections $S$ (which we will also denote by $S$), whose elements all have order 2. Measuring the distance between two chambers $gB$ and $hB$ boils down to finding an element $\tilde{w}\in W$ such that $g^{-1}h = bwb'$, for $b,b'\in B$, something we will not prove.

The group $G$ admits a \emph{Bruhat decomposition}. This means that $G$ decomposes as 
\[G = \bigsqcup_{\tilde{w}\in W}BwB,\]
where $w$ is a choice of representative of $\tilde{w}$. That this covers $G$ and that the different $BwB$ are disjoint is far from obvious. We refer the reader to \cite[\S6.1]{AbramBrown}. This decomposition implies that $(gB)^{-1}hB = Bg^{-1}hB = BwB$ for some $w\in N$. We claim without proof that $gB$ and $hB$ are adjacent if and only if $w = s$ for a simple reflection $s\in S\subset N$. This also implies the criterion for larger distances: $gB$ and $hB$ are a distance $r$ apart, if and only if $g^{-1}h\in Bs_{i_1}\cdots s_{i_r}B$ with $s_{i_j}\in S$. The decomposition also immediately tells us that all cosets of $B$ take the form
\[bwB,\ b\in B, w\in N.\]

\subsection{The representatives of $ G / B $}
We now give an overview of the method we developed for finding canonical representatives for the cosets $gB$. The novelty in the method is the use of the \emph{filtration matrix} to easily read off the structure of certain conjugates of the group $B$. It offers a method for finding representatives of $G/B$ that is very light on calculations.

Interestingly enough -- but not at all coincidentally -- the description of $G / B$ for an affine building $\mc B (SL_n(\Q),\nu_p)$ gives us a description of the chamber set for a spherical building $\mc B_{sph}( GL_n(\Z/p\Z))$ over a finite field for free. Put simply, we can perform our method below, but take coordinates modulo $p\Z$ in the necessary places. We will therefore take the building $\mc B (SL_3 (\Q),\nu_p)$ as the example that we work out in more detail here, but hold in mind that the method generalises to other affine buildings.
\bigskip

\noindent From the previous section it is clear that $B$ acts transitively on the cosets contained in $BwB$ by left multiplication. We fix a $w\in N$ and we try to find elements of $B$ that parametrise these cosets. The stabiliser of both $B$ and $wB$ is equal to $B\cap wBw^{-1}$. By the orbit stabiliser theorem the orbit of $wB$ under the action of $B$ is in bijection with $B/\left( B\cap wBw^{-1}\right)$. This is an important fact that we will use to determine a set of representatives of $BwB/B$.

We recall the definition of the subgroup $B_p$ of $SL_3(\Q)$:
\[B_p = \left\{\begin{pmatrix}  \nu_p = 0 & \nu_p \geq 0 & \nu_p \geq 0\\
                                \nu_p > 0 & \nu_p = 0 & \nu_p \geq 0\\
                                \nu_p > 0 & \nu_p > 0 & \nu_p = 0\end{pmatrix}\right\}.\]
It has an interesting, and rather simple, factorisation. We define a set of six matrix groups $U_{(i,j),r}$ with $0\leq i,j\leq 2$, $i\not=j$ and $r\in \Z\cup\{\infty\}$. An element $u\in U_{(i,j),r}$ has only ones on the diagonal, an element $q\in \Q$ in position $(i,j)$ such that $\nu_p(q) \geq r$ and zeros in all other positions. Examples of elements of $U_{(0,1),3}$, if $p\not= 2$, are the matrices
\[\begin{pmatrix}   1 & p^3 & 0\\
                    0 & 1 & 0\\
                    0 & 0 & 1\end{pmatrix},\ 
  \begin{pmatrix}   1 & p^3/2 & 0\\
                    0 & 1 & 0\\
                    0 & 0 & 1\end{pmatrix}\text{ and } 
  \begin{pmatrix}   1 & p^5 & 0\\
                    0 & 1 & 0\\
                    0 & 0 & 1\end{pmatrix}.\]
It turns out that $B_p$ admits a factorisation into these groups and the subgroup $T = N\cap B_p$, as follows. Define an element $F_{B_p}$ of $B_p$ by
\[F_{B_p} = \begin{pmatrix} 1 & 1 & 1\\
                        p & p-1 & 1\\
                        p & p & p-1\end{pmatrix}\]
which we call the \emph{filtration matrix of $B_p$}. The point of this matrix is that it is a choice of an element of $B_p$ whose coordinates have the minimal allowed $p$-adic valuation in $B_p$. It comes with a function
\[\phi_{B_p}(i,j)=\nu_p ( (F_B )_{i,j}),\]
i.e. $\phi_{B_p} ( i,j)$ is the valuation of the element in position $(i,j)$ of $F_{B_p}$. The multiplication map
\[T\times\prod_{i\not=j} U_{(i,j),\phi_{B_p}(i,j)}\to B_p\]
is bijective, although not a group homomorphism.

Since we are looking for representatives of $B_p/ (B_p\cap wB_pw^{-1})$, we want to have a similar description for $wB_pw^{-1}$. The definition of a filtration matrix is easily extended to all $wB_pw^{-1}$, simply by
\[F_{wB_p} = wF_{B_p}w^{-1},\]
where we use $wB_p$ in the subscript since $wB_pw^{-1}$ is the stabiliser of the chamber $wB_p$. This too comes with a map $\phi_{wB_p}(i,j) = \nu_p ( (F_{wB_p})_{i,j})$. Since the simple reflections only permute rows and columns and multiply them by some scalar, but do not add rows or columns together, $F_{wB_p}$ is again an element of $wB_pw^{-1}$ whose coordinates have minimal valuation. Therefore, $wB_pw^{-1}$ is characterised by $\phi_{wB_pw^{-1}}$, and we again have a bijection
\[T\times\prod_{i\not=j} U_{(i,j),\phi_{wB_p}(i,j)}\to wB_pw^{-1}.\]
Note that this also implies that
\[B_p\cap wB_pw^{-1} \cong T\times\prod_{i\not=j} U_{(i,j),r^w_{i,j}}\]
with $r^w_{i,j} = \max(\phi_{wB_p}(i,j),\phi_{B_p}(i,j))$.

The groups $U_{(i,j),r}$ are isomorphic to subgroups of $\Q$. Namely, we can send a matrix in $U_{(i,j),r}$ to the element in $\Q$ corresponding to the coordinate in position $(i,j)$. By this map, for every $r$ the group $U_{(i,j),r}$ is isomorphic to $p^r\Z_{(p)}$ as an abelian group, where
\[\Z_{(p)} = \{z\in \Q\ |\ \nu_p (z) \geq 0\}.\]
So $p^r\Z_{(p)}$ is exactly the abelian group of elements with $p$-adic valuation greater then or equal to $r$. A quotient $p^r\Z_{(p)}/p^{r+n}\Z_{(p)}$ is isomorphic to $\Z/p^n\Z$.

Calculating representatives of the quotient $B / ( B\cap wBw^{-1})$ is an easy task now:
\[\begin{split}
    B / (B\cap wBw^{-1})    &\cong (T\times\prod_{i\not=j} U_{(i,j),\phi_{B_p}(i,j)}) / (T\times\prod_{i\not=j} U_{(i,j),r^w_{i,j}})\\
                            &\cong \prod_{i\not=j} (U_{(i,j),\phi_{B_p}(i,j)}/U_{(i,j),r^w_{i,j}})\\
                            &\cong \prod_{i\not=j} p^{\phi_{B_p} (i,j)}\Z_{(p)}/ p^{r_{i,j}^w}\Z_{(p)}\\
                            &\cong \prod_{i\not=j} \Z / p^{r_{i,j}^w - \phi_{B_p} (i,j)}\Z,
\end{split}\]
where the first isomorphism is an isomorphism of sets, and all others are isomorphisms of groups.
\bigskip

\noindent All of this leads to the conclusion that we can algorithmically parametrise the chambers in $BwB$ by products of elements of the groups
\begin{equation}\label{eq:rootquotient}
    U_{(i,j),\phi_{B_p}(i,j)}/U_{(i,j),r_{i,j}^w}\cong \Z/p^{r^w_{i,j}-\phi_{B_p}(i,j)}\Z.\tag{$\ast$}
\end{equation}
If $\phi_{B_p}(i,j) = 0$, so if $i<j$, we take a generator of such a group to be a matrix with only ones on the diagonal and at the coordinate $(i,j)$, and zeros everywhere else. If $\phi_{B_p}(i,j) = 1$, so if $i>j$, we take it to be the matrix with ones on the diagonal, $p$ in position $(i,j)$, and zeros everywhere else. The generator has order $p^{r_{i,j}^w}$ if $i<j$ in \eqref{eq:rootquotient} and order $p^{r_{i,j}^w-1}$ otherwise.

\begin{figure}[hbt]
    \centering
    \subfloat[A small part of the building $\mc B (SL_3(\Q),\nu_2)$ up to distance 1 from the chamber $B$. Pairs of chambers indicated by $Bs_iB$ are the cosets of $B$ in $Bs_iB$.]{
    \includegraphics[width = 0.55\textwidth]{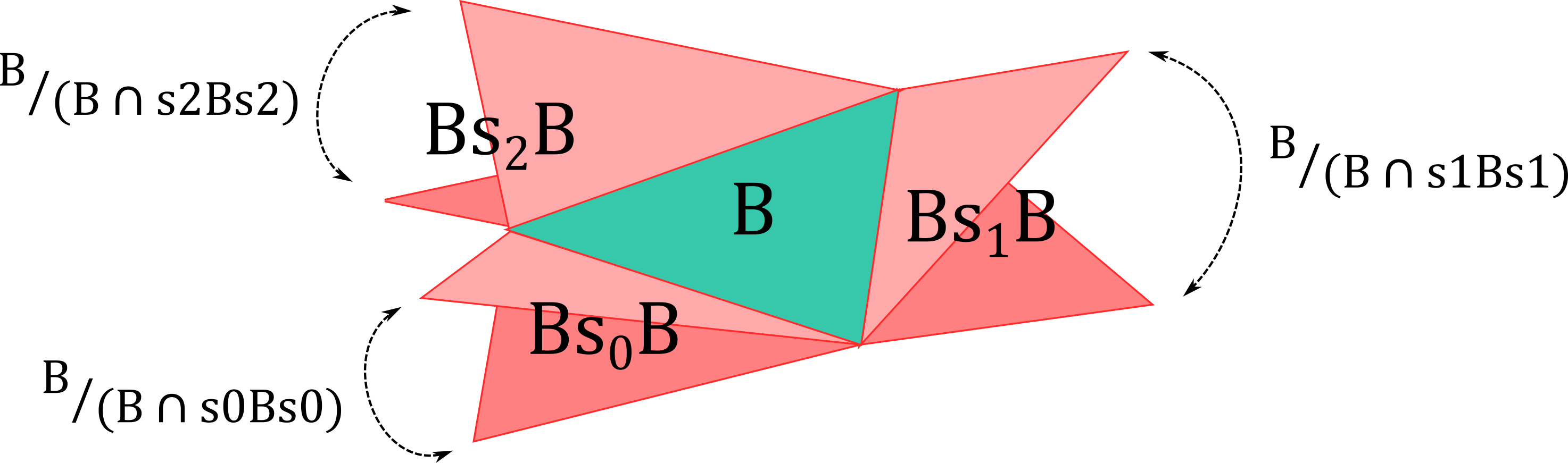}
    \label{fig:BuildAffA2a}}\qquad
    \subfloat[The same part of the building $\mc B (SL_3(\Q),\nu_2)$ with the labels of the chambers in place.]{
    \includegraphics[width = 0.35\textwidth]{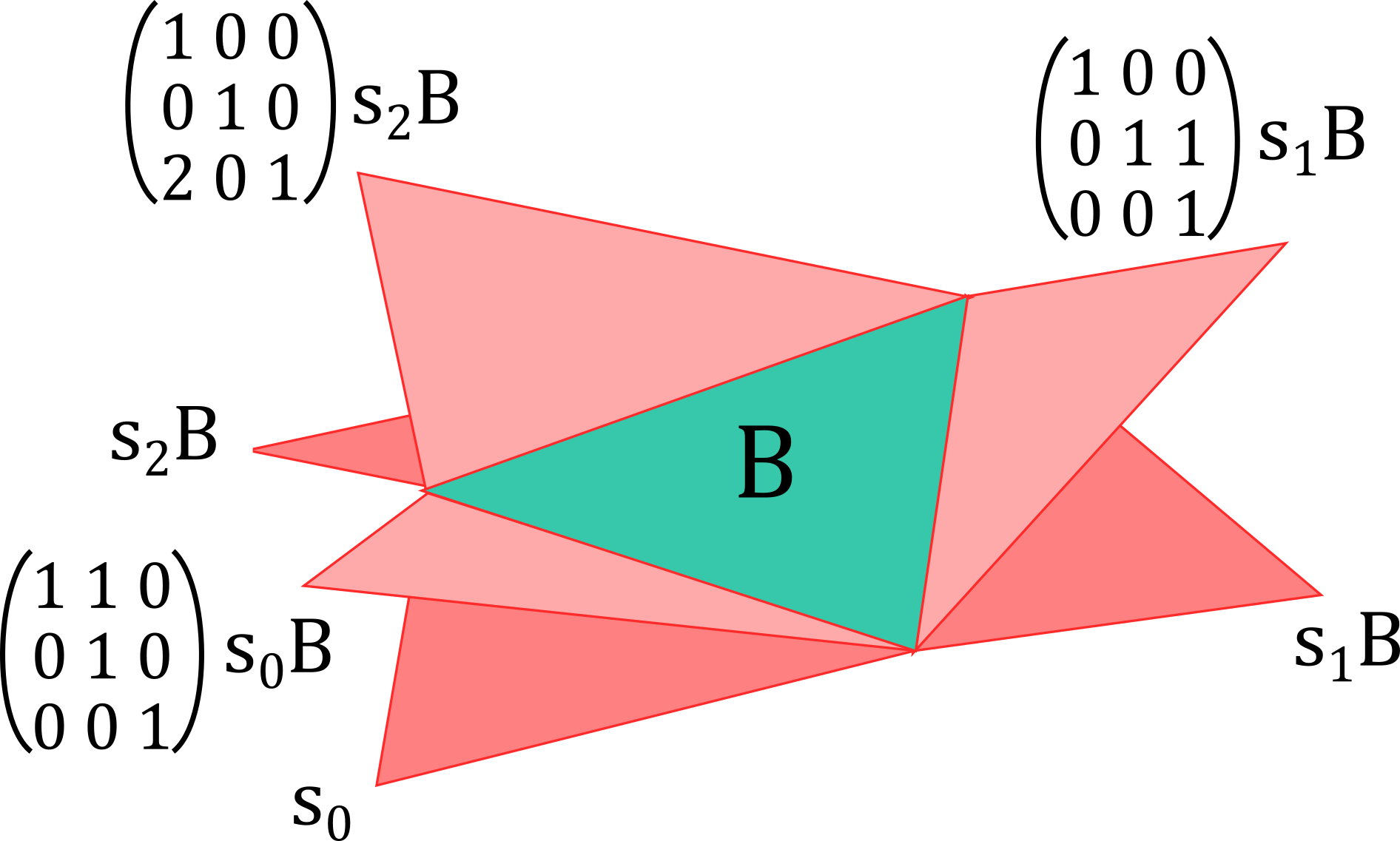}
    \label{fig:BuildAffA2b}}
    \caption{}
    \label{fig:BuildAffA2}
\end{figure}

Figures \ref{fig:BuildAffA2}\subref{fig:BuildAffA2a} and  \subref{fig:BuildAffA2b} show what this looks like in our example for the chambers closest to $B$.

%\bibliographystyle{apacite}
%\bibliography{./Biblio}

\begin{thebibliography}{}

\bibitem [\protect \citeauthoryear {%
Abramenko%
\ \BBA {} Brown%
}{%
Abramenko%
\ \BBA {} Brown%
}{%
{\protect \APACyear {2008}}%
}]{%
AbramBrown}
\APACinsertmetastar {%
AbramBrown}%
\begin{APACrefauthors}%
Abramenko, P.%
\BCBT {}\ \BBA {} Brown, K\BPBI S.%
\end{APACrefauthors}%
\unskip\
\newblock
\APACrefYear{2008}.
\newblock
\APACrefbtitle {Buildings, Theory and Applications} {Buildings, theory and
  applications}.
\newblock
\APACaddressPublisher{}{Springer, New-York}.
\PrintBackRefs{\CurrentBib}

\bibitem [\protect \citeauthoryear {%
Bekker%
}{%
Bekker%
}{%
{\protect \APACyear {2021}}%
}]{%
Bekker2020}
\APACinsertmetastar {%
Bekker2020}%
\begin{APACrefauthors}%
Bekker, B.%
\end{APACrefauthors}%
\unskip\
\newblock
\APACrefYearMonthDay{2021}{Jun}{}.
\newblock
\APACrefbtitle {The {B}uildings {G}allery.} {The {B}uildings {G}allery.}
\newblock
\begin{APACrefURL} [{2021-06-29}]\url{https://buildings.gallery/}
  \end{APACrefURL}
\PrintBackRefs{\CurrentBib}

\bibitem [\protect \citeauthoryear {%
Bruhat%
\ \BBA {} Tits%
}{%
Bruhat%
\ \BBA {} Tits%
}{%
{\protect \APACyear {1972}}%
}]{%
BrTi1}
\APACinsertmetastar {%
BrTi1}%
\begin{APACrefauthors}%
Bruhat, F.%
\BCBT {}\ \BBA {} Tits, J.%
\end{APACrefauthors}%
\unskip\
\newblock
\APACrefYearMonthDay{1972}{}{}.
\newblock
{\BBOQ}\APACrefatitle {Groupes r\'{e}ductifs sur un corps local} {Groupes
  r\'{e}ductifs sur un corps local}.{\BBCQ}
\newblock
\APACjournalVolNumPages{Inst. Hautes \'{E}tudes Sci. Publ.
  Math.}{}{41}{5--251}.
\PrintBackRefs{\CurrentBib}

\bibitem [\protect \citeauthoryear {%
Buekenhoudt%
}{%
Buekenhoudt%
}{%
{\protect \APACyear {2014}}%
}]{%
TitsBio2014}
\APACinsertmetastar {%
TitsBio2014}%
\begin{APACrefauthors}%
Buekenhoudt, F.%
\end{APACrefauthors}%
\unskip\
\newblock
\APACrefYearMonthDay{2014}{}{}.
\newblock
{\BBOQ}\APACrefatitle {A {B}iography of {J}acques {T}its} {A {B}iography of
  {J}acques {T}its}.{\BBCQ}
\newblock
\BIn{} R.~Holden Helge;~Piene\ (\BED), \APACrefbtitle {The {A}bel {P}rize
  2008-2012} {The {A}bel {P}rize 2008-2012}\ (\BPG~35-53).
\newblock
\APACaddressPublisher{}{Springer-Verlag}.
\PrintBackRefs{\CurrentBib}

\bibitem [\protect \citeauthoryear {%
Garrett%
}{%
Garrett%
}{%
{\protect \APACyear {1997}}%
}]{%
Garrett1997}
\APACinsertmetastar {%
Garrett1997}%
\begin{APACrefauthors}%
Garrett, P.%
\end{APACrefauthors}%
\unskip\
\newblock
\APACrefYear{1997}.
\newblock
\APACrefbtitle {Buildings and Classical Groups} {Buildings and classical
  groups}\ (\PrintOrdinal{1}\ \BEd).
\newblock
\APACaddressPublisher{}{Springer}.
\PrintBackRefs{\CurrentBib}

\bibitem [\protect \citeauthoryear {%
Schneider%
\ \BBA {} Stuhler%
}{%
Schneider%
\ \BBA {} Stuhler%
}{%
{\protect \APACyear {1997}}%
}]{%
SchneiStuhl}
\APACinsertmetastar {%
SchneiStuhl}%
\begin{APACrefauthors}%
Schneider, P.%
\BCBT {}\ \BBA {} Stuhler, U.%
\end{APACrefauthors}%
\unskip\
\newblock
\APACrefYearMonthDay{1997}{}{}.
\newblock
{\BBOQ}\APACrefatitle {Representation theory and sheaves on the Bruhat-Tits
  building} {Representation theory and sheaves on the bruhat-tits
  building}.{\BBCQ}
\newblock
\APACjournalVolNumPages{Publications Math\'ematiques de
  l'IH\'ES}{85}{}{97-191}.
\newblock
\begin{APACrefURL} \url{http://www.numdam.org/item/PMIHES\_1997\_\_85\_\_97\_0}
  \end{APACrefURL}
\PrintBackRefs{\CurrentBib}

\bibitem [\protect \citeauthoryear {%
Tits%
}{%
Tits%
}{%
{\protect \APACyear {1997}}%
}]{%
Tits1997}
\APACinsertmetastar {%
Tits1997}%
\begin{APACrefauthors}%
Tits, J.%
\end{APACrefauthors}%
\unskip\
\newblock
\APACrefYearMonthDay{1997}{}{}.
\newblock
{\BBOQ}\APACrefatitle {Reductive groups over local fields} {Reductive groups
  over local fields}.{\BBCQ}
\newblock
\APACjournalVolNumPages{Proceedings of Symposia in Pure
  Mathematics}{33}{}{29-69}.
\PrintBackRefs{\CurrentBib}

\end{thebibliography}

\end{document}